\documentclass[12pt]{article}
\usepackage{amssymb,amsmath}
\usepackage{cases}
\usepackage{amsfonts}
\usepackage{color,xcolor}
\usepackage[left=2.0cm,right=2.0cm,top=2.0cm,bottom=2.0cm]{geometry}
\usepackage[colorlinks,citecolor=blue,urlcolor=blue]{hyperref}

\newtheorem{theorem}{Theorem}[section]

\newtheorem{lemma}{Lemma}[section]

\newtheorem{definition}{Definition}[section]
\newtheorem{remark}{Remark}[section]
\newtheorem{example}{Example}[section]

\newcommand{\bal}{\begin{align}}
\newcommand{\bbal}{\begin{align*}}
\newcommand{\beq}{\begin{equation}}
\newcommand{\eeq}{\end{equation}}
\newcommand{\bca}{\begin{cases}}
\newcommand{\eca}{\end{cases}}
\newcommand{\pa}{\partial}
\newcommand{\fr}{\frac}
\newcommand{\na}{\nabla}
\newcommand{\De}{\Delta}

\newcommand{\cd}{\cdot}
\newcommand{\ep}{\varepsilon}
\newcommand{\dd}{\mathrm{d}}
\newcommand{\B}{\dot{B}}

\newcommand{\R}{\mathbb{R}}

\newcommand{\les}{\lesssim}

\newcommand{\D}{\mathrm{div}}
\newcommand{\n}{\nabla}
\newcommand{\bi}{\Big}

\begin{document}
\title{Remarks on the global large solution to the three-dimensional incompressible Navier-Stokes equations}

\author{Jinlu Li$^{1}$~ Yanghai Yu$^{2}$~and ~ Zhaoyang Yin$^{3,4}$ \thanks{E-mail: lijinlu@gnnu.cn (J. Li); yuyanghai214@sina.com (Y. Yu); mcsyzy@mail.sysu.edu.cn (Z. Yin)}\\
\small $^1$\it School of Mathematics and Computer Sciences, Gannan Normal University, Ganzhou 341000, China\\
\small $^2$\it School of Mathematics and Statistics, Anhui Normal University, Wuhu, Anhui, 241002, China\\
\small $^3$\it Department of Mathematics, Sun Yat-sen University, Guangzhou 510275, China\\
\small $^4$\it Faculty of Information Technology, Macau University of Science and Technology, Macau, China\\}

%\date{\today}

\maketitle\noindent{\hrulefill}

{\bf Abstract:}  In this paper, we derive a new smallness hypothesis of initial data for the three-dimensional incompressible Navier-Stokes equations. That is, we prove that there exist two positive constants
$c_0,C_0$ such that if
\begin{equation*}
  \|u_0^1+u^2_0,u^3_0\|_{\dot{B}_{p,1}^{-1+\frac{3}{p}}}  \|u^1_0,u^2_0\|_{\dot{B}_{p,1}^{-1+\frac{3}{p}}} \exp\{C_0 (\|u_0\|^{2}_{\dot{B}_{\infty,2}^{-1}}+\|u_0\|_{\dot{B}_{\infty,1}^{-1}})\} \leq c_0,
\end{equation*}
then \eqref{NS} has a unique global solution. As an application we construct two family of smooth solutions to the Navier-Stokes equations whose $\B^{-1}_{\infty,\infty}(\mathbb{R}^3)$ norm can be arbitrarily large.


\vskip0mm\noindent{\hrulefill}

\section{Introduction}\label{sec1}
In this paper, we focus on the incompressible Navier-Stokes equations in the whole space $\R^3$
\begin{eqnarray}\label{NS}
        (\rm{NS})~~\left\{\begin{array}{ll}
          \partial_tu+u\cd\na u-\De u+\na p=0,\\
          \D u=0,\\
          u(0,x)=u_0.
          \end{array}\right.
        \end{eqnarray}
where $u=(u^1(t,x),u^2(t,x),u^3(t,x))\in\R^3$ denotes the divergence free velocity field and $p\in \R$ is the scalar pressure.

Firstly, it should be noticed that the incompressible Navier-Stokes system \eqref{NS} is translation and scaling invariant: if $(u,\pi)$ is a solution of \eqref{NS} on $[0,T]\times\R^3$, then, for any positive $\lambda$, the scaled triplet $(u,\pi)_\lambda$ defined by
\bal\label{s}
(u,\pi)_\lambda(x,t)=(\lambda u(\lambda x,\lambda^2 t),\lambda^2\pi(\lambda x,\lambda^2 t))
\end{align}
is also a solution of \eqref{NS} on $[0,\lambda^{-2}T]\times\R^3$. Thus, this leads to the notation of critical spaces for well-posedness, that is, an adapted space must be translation and scaling invariant in the following sense: $||(u,\pi)_\lambda||_{X}\thickapprox||(u,\pi)||_{X}.$
The reader may check that the following spaces have the above invariance
$$\dot{H}^{\fr12}(\R^3)\hookrightarrow L^3(\R^3)\hookrightarrow\dot{B}_{p,1}^{\frac{3}{p}-1}(\R^3)|_{p<\infty}\hookrightarrow BMO^{-1}(\R^3)\hookrightarrow\B^{-1}_{\infty,\infty}(\R^3).$$

Now, let us recall some important progress about the global existence results for small data. In his seminal work, Leray \cite{Leray} proved in 1934 that if $\|u_0\|_{L^2}\|\n u_0\|_{L^2}$ is small enough, then there exists a global regular solution of (NS). Then Fujita--Kato \cite{Fujita1964} proved in 1964 that if $\|u_0\|_{\dot{H}^\fr12}$ is small enough, then there exists a unique global solution in the space $C_b(\R^+;\dot{H}^\fr12)\cap L^4(\R^+;\dot{H}^1)$. Indeed, the theorem of Fujita--Kato \cite{Fujita1964} allows to construct local in time unique solution to \eqref{NS} with initial data in the Lebsegue space $L^3(\R^3)$ \cite{Kato}. Moreover, if the initial norm $\|u_0\|_{L^3}$ is sufficiently small, then the strong solution to \eqref{NS} exists globally in time. The above result was extended by Cannone--Meyer--Planchon \cite{Cannone 1994} for initial data in
Besov spaces with negative index. More precisely, they proved that if the initial data belongs
to the Besov space, $\dot{B}_{p,1}^{\frac{3}{p}-1}(\R^3)$ for some $p \in(3,\infty)$ and its norm is sufficiently small, then \eqref{NS} has a unique global solution. The typical example of such kind
of initial data reads
\bbal
u_{0,\ep}(x)=\ep^{-\alpha}\sin\bi(\frac{x_3}{\ep}\bi)(\pa_2\phi,\pa_1\phi,0)(x)\quad\mbox{with}\quad\alpha\in(0,1),\quad\phi\in \mathcal{S}(\R^3;\R).
\end{align*}
We remark that this type of initial data is not small in either $\dot{H}^{\fr12}(\R^3)$ or $L^3(\R^3)$.

Koch--Tataru in 2001 \cite{Koch2001} proved that given initial data in the derivatives of BMO space and its norm $BMO^{-1}$ is sufficiently small, then \eqref{NS} has a unique global solution. We point out that the largest space and the
norm of which is scaling invariant under \eqref{s}, is $\B^{-1}_{\infty,\infty}(\R^3)$. Hmidi--Li \cite{Hmidi} showed that smallness of $\dot{B}^{-1}_{\infty,\infty}$ norm of solution to d-dimensional ($d\geq3$) incompressible Navier-Stokes prevents blowups. Moreover, Bourgain--Pavlovi\'{c}
\cite{Bourgain} proved that \eqref{NS} is actually ill-posed with initial data in $\B^{-1}_{\infty,\infty}(\R^3)$. Let us simply notice that the above norms have the following relation
\bbal
\|u_0\|_{\B^{-1}_{\infty,\infty}(\R^3)}\lesssim\|u_0\|_{BMO^{-1}(\R^3)}\lesssim\|u_0\|_{\B^{-1}_{\infty,2}(\R^3)}
\end{align*}
with
\bbal
\|u_0\|_{\B^{-1}_{\infty,\infty}(\R^3)}=\sup_{t>0}t^\fr12||e^{t\Delta}u_0||_{L^\infty}\quad\mbox{and}\quad\|u_0\|_{\B^{-1}_{\infty,2}(\R^3)}=\|e^{t\Delta}u_0\|_{L^{2}(\R^+,L^\infty(\R^3))}
\end{align*}
 A remark due to Y. Meyer \cite{Meyer 1997} is that the norm in such a space is always greater than the norm in the Besov space $\B^{-1}_{\infty,\infty}(\R^3)$. This leads to the definition of a large initial data for the incompressible Navier-Stokes equations: A divergence free vector field $u_0$ is a large initial data for the incomcpressible Navier-Stokes system if its $\B^{-1}_{\infty,\infty}$ norm is large. Chemin--Gallagher \cite{Chemin2009} proved that if a certain nonlinear function of the initial data is small enough, then there is a global solution to the Navier-Stokes equations \eqref{NS} in a Koch-Tataru type space. Meanwhile, they provide an example of initial data satisfying that nonlinear smallness condition, but whose norm is arbitrarily large in $\B^{-1}_{\infty,\infty}(\R^3)$. For more results of large initial data which generate unique global solutions to \eqref{NS}, we refer the reader to see \cite{Chemin2006,Chemin2010,Chemin2013,Chemin2015,Lei2015,Liu2018} and the references therein. Recently, Li--Yu--Zhu--Yin \cite{Li2019} obtained the global large solution for a special initial data with the first two component of the initial velocity field being large in Besov space $\B^{-1}_{\infty,\infty}(\R^3)$. One can mention that Paicu--Zhang \cite{Zhang2014} proved the
global well-posedness of incompressible inhomogeneous Navier-Stokes equation with the third component of the initial velocity field being
large. Motivated by this work \cite{Chemin2009,Zhang2014,Li2019}, we continue to study the global behavior of solutions to the incompressible Navier-Stokes equations \eqref{NS} for arbitrary large initial in the present paper.

The main result of the paper read as follows:
\begin{theorem}\label{the1}
Let $3<p<6$.  Then there exist two constants $\delta,C>0$ such that for any $u_0=(u_0^1,u^2_0,u_0^3) \in \dot{B}_{p,1}^{\frac{3}{p}-1}$ satisfying the condition
\begin{equation}\label{con}
  \|u_0^1+u^2_0,u^3_0\|_{\dot{B}_{p,1}^{-1+\frac{3}{p}}}  \|u^1_0,u^2_0\|_{\dot{B}_{p,1}^{-1+\frac{3}{p}}} \exp\{C (\|u_0\|^{2}_{\dot{B}_{\infty,2}^{-1}}+\|u_0\|_{\dot{B}_{\infty,1}^{-1}})\} \leq \delta,
\end{equation}
then \eqref{NS} admits a unique global solution $u$.
\end{theorem}

\begin{remark}
Our obtained result in Theorem 1.1 improves considerably the corresponding result in \cite{Zhang2010} when the initial data belongs to Besov spaces $B^{-1+\frac 3p}_{p,1},3<p<6$. From Example 1.1, we can obtain the global solution of the initial data as \eqref{initial data1}. This implies that our obtained result in Theorem 1.1 partially covers the result in \cite{Chemin2009} when $\alpha\in(\frac67,1)$, and our proof is more brief than that in \cite{Chemin2009}.
\end{remark}

Let us present some examples of initial data the norm of which are big in $\B^{-1}_{\infty,\infty}(\R^3)$, yet they satisfy the smallness condition \eqref{con}.
The first case is that the third component of the initial velocity field is large, and the first two components of the initial velocity field are small.

\begin{example}
Let $p\in [5,6)$ and $\alpha\in(\frac{6}{p+2},1)$. According to \cite{Chemin2009}, we take the following initial data
\bal\label{initial data1}
u_{0,\ep}(x)=\bi(\log\frac1\ep\bi)^{\frac15}\ep^{-1}\cos\bi(\frac{x_1}{\ep}\bi)(0,-\ep^{\alpha}\pa_3\phi,\pa_2\phi)(x_1,x_2/\ep^{\alpha},x_3).
\end{align}
From Lemma 3.1 of \cite{Chemin2009}, we have
\bbal
&||\cos\bi(\frac{x_1}{\ep}\bi)\na\phi(x_1,x_2/\ep^{\alpha},x_3)||_{\B^{-1+\frac3p}_{p,1}}\lesssim \ep^{1-\frac3p+\frac\alpha p},
\\& ||\cos\bi(\frac{x_1}{\ep}\bi)\na\phi(x_1,x_2/\ep^{\alpha},x_3)||_{\B^{-1}_{\infty,\infty}}\gtrsim\ep.
\end{align*}
Thus, we can get
\bbal
&||u_{0,\ep}||_{\B^{-1}_{\infty,\infty}}\approx ||u_{0,\ep}||_{\B^{-1}_{\infty,1}}\approx \bi(\log\frac1\ep\bi)^{\frac15},
\\&||u^3_{0,\ep}||_{\B^{-1+\frac3p}_{p,1}}\lesssim \ep^{-\frac3p+\frac\alpha p}\bi(\log\frac1\ep\bi)^{\frac15},
\\&||u^1_{0,\ep},u^2_{0,\ep}||_{\B^{-1+\frac3p}_{p,1}}\lesssim \ep^{\alpha-\frac3p+\frac\alpha p}\bi(\log\frac1\ep\bi)^{\frac15}.
\end{align*}
Then, direct calculations show that the left side of \eqref{con} becomes
\begin{equation*}
C\bi(\log\frac1\ep\bi)^{\frac25}.\exp\bi(C(\log \frac1\ep)^{\frac25}\bi)\ep^{\alpha-\frac6p+\frac{2\alpha}{p}},
\end{equation*}
which implies \eqref{NS} have a global solution  for $\ep$ sufficiently small.
\end{example}

The second case is that the first two components of the initial velocity field are large and the third component of the initial velocity field is small.
\begin{example}
Let $p\in(3,6)$. According to \cite{Li2019}, we take the initial data $u_{0,\ep}=(\pa_2a_{0,\ep},-\pa_1a_{0,\ep},0)$ with
\bbal
a_{0,\ep}(x_1,x_2,x_3)=\ep^{-1}\bi(\log\log\frac1\ep\bi)^{\frac12}  \chi(x_1,x_2)\phi(x_3),
\end{align*}
where the smooth functions $\chi,\phi$ satisfying $\hat{\chi}(-\xi_1,-\xi_2)=\hat{\chi}(\xi_1,\xi_2)$, $\hat{\phi}(-\xi_3)=\hat{\phi}(\xi_3)$,
\begin{align*}
\mathrm{supp} \hat{\chi}\subset \mathcal{\widetilde{C}},\quad \hat{\chi}(\xi_1,\xi_2)\in[0,1]; \quad \hat{\chi}(\xi_1,\xi_2)=1 \quad\mbox{for} \quad (\xi_1,\xi_2)\in\mathcal{\widetilde{C}}_1,
\end{align*}
and
\begin{align*}
\hat{\phi}(\xi_3)=0\quad \mbox{for} \quad|\xi_3|\in\bi[\frac{2\sqrt{2}}{3},\frac{3\sqrt{2}}{4}\bi]^c,\quad \hat{\phi}(\xi_3)\in[0,1]; \quad \hat{\phi}(\xi_3)=1 \quad\mbox{for}  \quad |\xi_3|\in\bi[\frac{\sqrt{34}}{6},\frac{\sqrt{17}}{4}\bi],
\end{align*}
where
\begin{align*}
&\mathcal{\widetilde{C}}\triangleq\Big\{\xi\in\R^2: \ |\xi_1-\xi_2|\leq \ep,\ \frac89\leq\xi^2_1+\xi^2_2\leq \frac98\Big\},
\\&\mathcal{\widetilde{C}}_1\triangleq\Big\{\xi\in\R^2: \ |\xi_1-\xi_2|\leq \frac\ep2,\ \frac{17}{18}\leq\xi^2_1+\xi^2_2\leq \frac{17}{16}\Big\} .
\end{align*}
In fact, one has
\begin{align*}
||\hat{a}_0||_{L^{\frac{p}{p-1}}}\approx \ep^{-\frac{1}{p}}\Big(\log\log\frac1\ep\Big)^\frac12.
\end{align*}
Then, direct calculations show that the left side of \eqref{con} becomes
\begin{align*}
C\ep^{1-\frac2p}\bi(\log\log \frac1\ep\bi)\exp\bi(C\log\log \frac1\ep\bi),
\end{align*}
which implies \eqref{NS} have a global solution  for $\ep$ sufficiently small. From Remark 2.1 of  \cite{Li2019}, it also holds
\bbal
||u_0||_{\dot{B}^{-1}_{\infty,\infty}}\approx ||u_0||_{L^\infty}\gtrsim \bi(\log\log \frac1\ep\bi)^\frac12.
\end{align*}
\end{example}

\section{Littlewood-Paley Analysis}
Throughout this paper, we will denote by $C$ any constant which may change from line to line and write $A\lesssim B$ if $A\leq CB$. $A\approx B$ means that $A\lesssim B$ and $B\lesssim A$. We also shall use the abbreviated notation $||f_1,\cdots,f_n||_{X}=||f_1||_{X}+\cdots+||f_n||_{X}$ for some Banach space $X$.

Next, we recall the Littlewood-Paley theory, the definition of homogeneous Besov spaces and some useful properties.

Let us start by introducing the Littlewood-Paley decomposition. Choose a radial function $\varphi\in \mathcal{S}(\mathbb{R}^d)$ supported in ${\mathcal{C}}=\{\xi\in\mathbb{R}^d:\frac34\leq |\xi|\leq \frac83\}$ such that
\begin{align*}
\sum_{j\in \mathbb{Z}}\varphi(2^{-j}\xi)=1 \quad \mathrm{for} \ \mathrm{all} \ \xi\neq0.
\end{align*}
The frequency localization operator $\dot{\Delta}_j$ and $\dot{S}_j$ are defined by
\begin{align*}
\dot{\Delta}_jf=\varphi(2^{-j}D)f=\mathcal{F}^{-1}(\varphi(2^{-j}\cdot)\mathcal{F}f) \quad \mbox{and}\quad\dot{S}_jf=\sum_{k\leq j-1}\dot{\Delta}_kf \quad \mathrm{for} \quad j\in\mathbb{Z}.
\end{align*}
With a suitable choice of $\varphi$, one can easily verify that
\begin{align*}
\dot{\Delta}_j\dot{\Delta}_kf=0 \quad \mathrm{if} \quad |j-k|\geq2\quad \mbox{and}\quad\dot{\Delta}_j(\dot{S}_{k-1}f\dot{\Delta}_kf)=0 \quad  \mathrm{if} \quad  |j-k|\geq5.
\end{align*}
Next we recall Bony's decomposition from \cite{Bahouri2011}:
\begin{align*}
uv=\dot{T}_uv+\dot{T}_vu+\dot{R}(u,v),
\end{align*}
with
\begin{align*}
\dot{T}_uv=\sum_{j\in\mathbb{Z}}\dot{S}_{j-1}u\dot{\Delta}_jv, \quad \quad \dot{R}(u,v)=\sum_{j\in\mathbb{Z}}\dot{\Delta}_ju\widetilde{\Delta}_jv, \quad \quad \widetilde{\Delta}_jv=\sum_{|j'-j|\leq 1}\dot{\Delta}_{j'}v.
\end{align*}
\begin{definition}
We denote by $\mathcal{Z}'(\mathbb{R}^d)$ the dual space of $\mathcal{Z}(\mathbb{R}^d)$, where we set $$\mathcal{Z}(\mathbb{R}^d)=\bi\{f\in \mathcal{S}(\mathbb{R}^d):D^{\alpha}\hat{f}(0)=0;\ \forall \alpha\in \mathbb{N}^d\bi\}.$$
\end{definition}
Then we have the formal homogenous Littlewood-Paley decomposition $$f=\sum\limits_{j\in \mathbb{Z}}\dot{\Delta}_jf,\quad \forall f\in\mathcal{Z}'(\mathbb{R}^d).$$

The operators $\dot{\Delta}_j$ help us recall the definition of the homogenous Besov space (see \cite{Bahouri2011}).

\begin{definition}
Let $s\in \mathbb{R}$ and $1\leq p,r\leq \infty$. The homogeneous Besov space $\dot{B}^s_{p,r}$ is defined by
\begin{align*}
\dot{B}^s_{p,r}=\bi\{f\in \mathcal{Z}'(\mathbb{R}^d):||f||_{\dot{B}^s_{p,r}}\triangleq \Big|\Big|(2^{ks}||\dot{\Delta}_k f||_{L^p})_{k\in\mathbb{Z}}\Big|\Big|_{\ell^r}<+\infty\bi\}.
\end{align*}
\end{definition}
It should be noted that a distribution $f\in\dot{B}^s_{p,r}$ if and only if there exist a constant $C>0$ and a non-negative sequence $\{d_{k}\}_{k\in\mathbb{Z}}$ such that
\begin{eqnarray*}
\forall k\in\mathbb{Z},\quad \|\dot{\Delta}_kf\|_{L^p}\leq Cd_{k} 2^{-ks}\|f\|_{{\dot{B}_{p,r}^{s}}} \quad \mbox{with} \quad \|d_{k}\|_{\ell^r}=1.
\end{eqnarray*}

\begin{lemma}\label{lem2.1}
Let $t>0$ and $1\leq p,p_1,p_2,r,r_1,r_2\leq \infty$ satisfying
\bbal
\frac1p=\frac{1}{p_1}+\frac{1}{p_2},\qquad \frac1r=\frac{1}{r_1}+\frac{1}{r_2}, \qquad t=t_1+t_2.
\end{align*}
Then there exists a positive constant $C$ such that
\begin{align}
&||\dot{T}_fg||_{\B^s_{p,r}}\leq C||f||_{L^{p_1}}||g||_{\B^s_{p_2,r}},\label{01}\\
&||\dot{T}_fg||_{\B^{s-t}_{p,r}}\leq C||f||_{\B^{-t}_{p_1,r_1}}||g||_{\B^s_{p_2,r_2}},\label{02}
\\&||\dot{R}(f,g)||_{\B^t_{p,r}}\leq C||f||_{B^{t_1}_{p_1,r_1}}||g||_{B^{t_2}_{p_2,r_2}}.\label{03}
\end{align}
\end{lemma}
\textbf{Proof of Lemma \ref{lem2.1}}\quad The third conclusion \eqref{03} is the direct result of Theorem 2.52 in \cite{Bahouri2011}. The other conclusions come essentially from \cite{Bahouri2011}, we give the proof here for completeness. Then, by using the properties of spectral localization of the Littlewood-Paley decomposition, one has
\begin{align}\label{04}
    \|\dot{\Delta}_j{\dot{T}_u v}\|_{L^p}&\leq C\|\sum_{|k-j|\leq 4} \Delta_j(\dot{S}_{k-1}u \dot{\Delta}_k v)\|_{L^p} \nonumber \\
     &\leq C \sum_{|k-j|\leq 4} \|\dot{S}_{k-1}u\|_{L^{p_1}} \|\dot{\Delta}_k v\|_{L^{p_2}} \nonumber \\
          &\leq C 2^{-js}\|u\|_{L^{p_1}}\sum_{|k-j|\leq 4}   2^{(j-k)s}   2^{ks}\|\dot{\Delta}_k v\|_{L^{p_2}} \nonumber
     \\&\leq C d_j 2^{-j s} \|u\|_{L^{p_1}} \|v\|_{\dot{B}_{p_2,r}^{s}} \quad \mbox{with} \quad \|d_{j}\|_{\ell^r}=1
     \end{align}
and
\begin{align}\label{05}
    \|\dot{\Delta}_j{\dot{T}_u v}\|_{L^p}&\leq C\|\sum_{|k-j|\leq 4} \dot{\Delta}_j(\dot{S}_{k-1}u \dot{\Delta}_k v)\|_{L^p} \nonumber \\
     &\leq C \sum_{|k-j|\leq 4} \sum_{k' \leq k-2}  \|\dot{\Delta}_{k'}u\|_{L^{p_1}} \|\dot{\Delta}_k v\|_{L^{p_2}} \nonumber \\
          &\leq C 2^{-j(s-t)}\sum_{|k-j|\leq 4}  2^{(k-j)(-s+t)}\sum_{k' \leq k-2}  2^{(k'-k)t}  2^{-k't}\|\dot{\Delta}_{k'}u\|_{L^{p_1}} 2^{ks}\|\dot{\Delta}_k v\|_{L^{p_2}} \nonumber
     \\&\leq C d_j 2^{-j (s-t)} \|u\|_{\dot{B}_{p_1,r}^{-t}} \|v\|_{\dot{B}_{p_2,r}^{s}}.
     \end{align}
\eqref{04} and \eqref{05} result the desired \eqref{01} and \eqref{02}, respectively. This ends the proof of Lemma \ref{lem2.1}.
Next, we present the following product estimate which will be used in the sequel.
\begin{lemma}\cite{Danchin2001}\label{lem2.2}
Let $2\leq p\leq \infty$, $s_1\leq \frac dp$ and $s_2\leq \frac dp$ with $s_1+s_2>d\max\{0,\frac2p-1\}$. Then there holds
\begin{align*}
||fg||_{\dot{B}^{s_1+s_2-\frac dp}_{p,1}}&\leq C||f||_{\dot{B}^{s_1}_{p,1}}||g||_{\dot{B}^{s_2}_{p,1}}.
\end{align*}
\end{lemma}
Finally, we recall the optimal regularity estimates for the heat equations.
\begin{lemma}\cite{Danchin2005}\label{lem2.3}
Let $s\in \mathbb{R}$, $1\leq p,r\leq \infty$  and $1\leq q_1\leq q_2\leq \infty$. Assume that $u_0\in \dot{B}^s_{p,r}$ and $G\in {\widetilde{L}}^{q_1}_T(\dot{B}^{s+\frac{2}{q_1}-2}_{p,r})$. Then the heat equations
\begin{align*}
\left\{\begin{array}{ll}
\partial_tu-\Delta u=G,\\
 u(0,x)=u_0,
\end{array}\right.
\end{align*}
has a unique solution $u\in \widetilde{L}^{q_2}_T(\dot{B}^{s+\fr{2}{q_2}}_{p,r}))$ satisfying for all $T>0$
\begin{align*}
||u||_{\widetilde{L}^{q_2}_T(\dot{B}^{s+\fr{2}{q_2}}_{p,r})}\les||u_0||_{\dot{B}^s_{p,r}}+||G||_{{\widetilde{L}}^{q_1}_T(\dot{B}^{s+\frac{2}{q_1}-2}_{p,r})}.
\end{align*}
\end{lemma}

\section{Proof of the Main Result}\label{sec3}
\textbf{Proof of Theorem \ref{the1}.} Let $U=e^{t\Delta}u_0$ be the solutions generated by the following heat equations
\begin{eqnarray}\label{1jl}
        \left\{\begin{array}{ll}
          \pa_tU-\De U=0,\\
         \D U=0,\\
          U|_{t=0}=u_0.\end{array}\right.
\end{eqnarray}
Introducing the new quantity $v=u-U$, the system \eqref{NS} can be reduced to
\begin{eqnarray}\label{CNS}
        \left\{\begin{array}{ll}
          \partial_tv+v\cd\na v-\De v+\na p=-v\cd\na U-U\cd\na v-U\cd\na U,\\
          \D v=0, \\
           v_0=0.
          \end{array}\right.
\end{eqnarray}
To element the pressure term, applying the Leray operator $\mathbb{P}$ to the equation \eqref{NS}, one has
\begin{eqnarray}\label{3.1}
        \left\{\begin{array}{ll}
          \partial_tv-\De v=-\mathbb{P}(v\cd\na v-\D(v\otimes U)-\D(U\otimes v)-U\cd\na U),\\
          \D v=0, \\
           v_0=(0,0,0).
          \end{array}\right.
\end{eqnarray}
Invoking Lemma \ref{lem2.3} to the above system \eqref{3.1} yields
\bal\label{3.2}
||v||_{L^\infty_t(\B^{-1+\frac3p}_{p,1})}+||v||_{L^1_t(\B^{1+\frac3p}_{p,1})}
\les&~ \underbrace{\int^t_0||v\cd\na v||_{\B^{\frac3p-1}_{p,1}}\dd \tau}_{I_1}+\underbrace{\int^t_0||U\cd\na U||_{\B^{\frac3p-1}_{p,1}}\dd \tau}_{I_2}\nonumber\\
\quad&+\underbrace{\int^t_0||\D(v\otimes U)+\D(U\otimes v)||_{\B^{\frac3p-1}_{p,1}}\dd \tau}_{I_3},
\end{align}
where we have used the fact that $\mathbb{P}$ is a smooth homogeneous of degree 0 Fourier multipliers which maps $\dot{B}_{p,1}^{\frac3p-1}$ to itself.

For the term $I_1$, using the product estimate (see Lemma \ref{lem2.2}), we obtain
\bal\label{3.3}
I_1\les~\int_0^t||v||_{\B^{-1+\frac3p}_{p,1}}||v||_{\B^{1+\frac3p}_{p,1}}\dd \tau\les~||v||_{L^\infty_t(\B^{-1+\frac3p}_{p,1})}||v||_{L^1_t(\B^{1+\frac3p}_{p,1})}.
\end{align}
For the term $I_2$, notice that $\D U$=0, we have
\bbal
&U\cd\na U^1=(U^1+U^2)\pa_1U^1+U^2\pa_2(U^1+U^2)+U^2\pa_3U^3+U^3\pa_3U^1,
\\&U\cd\na U^2=(U^1+U^2)\pa_2U^2+U^1\pa_1(U^1+U^2)+U^1\pa_3U^3+U^3\pa_3U^2,
\\&U\cd\na U^3=U^1\pa_1U^3+U^2\pa_2U^3-U^3(\pa_1U^1+\pa_2U^2).
\end{align*}
Using the lemma \ref{lem2.2} again gives
\bal\label{3.4}
I_2=\int^t_0||U\cd\na U||_{\B^{-1+\frac3p}_{p,1}}\dd \tau\les&~ \int^t_0||U^1+U^2,U^3||_{\B^{\frac3p}_{p,1}}||U^1,U^2||_{\B^{\frac3p}_{p,1}}\dd\tau
\\\les&~ ||u_{0}^1+u_{0}^2,u_{0}^3||_{\B^{\frac3p-1}_{p,1}}||u_{0}^{1},u_{0}^{2}||_{\B^{\frac3p-1}_{p,1}}.
\end{align}
To deal with the term $I_3$, by Bony's decomposition, one has
\begin{equation*}
v\cd\na U^i=\sum^{3}_{j=1}\pa_j(U^iv^j)=\sum^{3}_{j=1}[\pa_j(\dot{T}_{U^i} v^j)+\dot{T}_{v^j} \pa_jU^i+\pa_j\dot{R}(U^i,v^j)]
\end{equation*}
and
\begin{equation*}
U\cd\na v^i=\sum^{3}_{j=1}\pa_j(U^jv^i)=\sum^{3}_{j=1}[\dot{T}_{U^j} \pa_jv^i+\pa_j(\dot{T}_{v^i} U^j)+\pa_j\dot{R}(U^j,v^i)].
\end{equation*}
Using \eqref{01}--\eqref{03} from Lemma \ref{lem2.1}, respectively, we obtain
\bal
&||\pa_j(\dot{T}_{U^i} v^j),\dot{T}_{U^j} \pa_jv^i||_{\B^{-1+\frac3p}_{p,1}}\les ~||U||_{L^\infty}||v||_{\B^{\frac3p}_{p,1}},\label{y1}\\
&||\dot{T}_{v^j} \pa_jU^i,\pa_j(\dot{T}_{v^i} U^j)||_{\B^{-1+\frac3p}_{p,1}}\les~||U||_{\B^1_{\infty,\infty}}||v||_{\B^{-1+\frac3p}_{p,1}},\label{y2}\\
&||\dot{R}(U^i,v^j),\dot{R}(U^j,v^i)||_{\B^{\frac3p}_{p,1}}\les~||U||_{\B^1_{\infty,\infty}}||v||_{\B^{-1+\frac3p}_{p,1}}.\label{y3}
\end{align}
Combining \eqref{y1}--\eqref{y3} implies
\bal\label{3.5}
\int^t_0||v\cd\na U+U\cd\na v||_{\B^{-1+\frac3p}_{p,1}}\dd \tau\les&~ \int^t_0||v||_{\B^{\frac3p}_{p,1}}||U||_{L^\infty}\dd \tau+\int^t_0||v||_{\B^{-1+\frac3p}_{p,1}}||U||_{\B^1_{\infty,\infty}}\dd \tau\nonumber\\
\les&~ \int^t_0||v||_{\B^{-1+\frac3p}_{p,1}}||U||_{\B^1_{\infty,\infty}}\dd \tau+\int^t_0||U||_{L^\infty}||v||^{\frac{1}{2}}_{\B^{-1+\frac3p}_{p,1}}||v||^{\frac{1}{2}}_{\B^{1+\frac3p}_{p,1}}\dd \tau\nonumber\\
\les&~ \int^t_0\bi(||U||^2_{L^\infty}+||U||_{\B^{1}_{\infty,\infty}}\bi)||v||_{\B^{-1+\frac3p}_{p,1}}\dd \tau+\frac12||v||_{L^1_t(\B^{1+\frac3p}_{p,1})}.
\end{align}
Putting the estimates \eqref{3.3}, \eqref{3.4} and \eqref{3.5} together with \eqref{3.2} yields
\bal\label{yyh}
||v||_{L^\infty_t(\B^{-1+\frac3p}_{p,1})}+||v||_{L^1_t(\B^{1+\frac3p}_{p,1})}
\les&~||v||_{L^\infty_t(\B^{-1+\frac3p}_{p,1})}||v||_{L^1_t(\B^{1+\frac3p}_{p,1})}+
||u_{0}^{1}+u_{0}^{2},u_{0}^{3}||_{\B^{\frac3p-1}_{p,1}}||u_{0}^{1},u_{0}^{2}||_{\B^{\frac3p-1}_{p,1}}
\nonumber\\
\quad& +\int^t_0\bi(||U||^2_{L^\infty}+||U||_{\B^{1}_{\infty,\infty}}\bi)||v||_{\B^{-1+\frac3p}_{p,1}}\dd \tau.
\end{align}
Now, we define
\bbal
\Gamma\triangleq\sup\bi\{t\in[0,T^*):||v||_{L^\infty_t(\B^{-1+\frac3p}_{p,1})}\leq \eta\ll1\bi\},
\end{align*}
where $\eta$ is a small enough positive constant which will be determined later on.

Thus, for all $t\in[0,\Gamma]$, choosing $\eta$ small enough, we infer from \eqref{yyh}
\bal\label{es-con1}
||v||_{L^\infty_t(\B^{-1+\frac3p}_{p,1})}+||v||_{L^1_t(\B^{1+\frac3p}_{p,1})}\les&~ ||u_{0}^{1}+u_{0}^{2},u_{0}^{3}||_{\B^{\frac3p-1}_{p,1}}||u_{0}^{1},u_{0}^{2}||_{\B^{\frac3p-1}_{p,1}}
\nonumber\\
\quad& +\int^t_0\bi(||U||^2_{L^\infty}+||U||_{\B^{1}_{\infty,\infty}}\bi)||v||_{\B^{-1+\frac3p}_{p,r}}\dd \tau.
\end{align}
Notice that $$\int^t_0\bi(||U||^2_{L^\infty}+||U||_{\B^{1}_{\infty,\infty}}\bi)\dd \tau\les ||u_0||^2_{\B^{-1}_{\infty,2}}+||u_0||_{\B^{-1}_{\infty,1}},$$
by Gronwall's inequality, we have for all $t\in[0,\Gamma]$
\bal\label{es-con2}
||v||_{L^\infty_t(\B^{\frac3p-1}_{p,1})}+||v||_{L^1_t(\B^{\frac3p+1}_{p,1})}\les&~||u_{0}^{1}+u_{0}^{2},u_{0}^{3}||_{\B^{\frac3p-1}_{p,1}}||u_{0}^{1},u_{0}^{2}||_{\B^{\frac3p-1}_{p,1}} \exp\bi\{C\bi(||u_0||^2_{\B^{-1}_{\infty,2}}+||u_0||_{\B^{-1}_{\infty,1}}\bi)\bi\}\nonumber\\
\leq&~C\delta
\end{align}
provided that the condition \eqref{con} holds.

Choosing $\eta=2C\delta$, thus we can get
\bbal
||v||_{L^\infty_t(\B^{-1+\frac3p}_{p,1})}&\leq \fr\eta2 \quad\mbox{for}\quad t\leq \Gamma.
\end{align*}

So if $\Gamma<T^*$, due to the continuity of the solutions, we can obtain that there exists $0<\epsilon\ll1$ such that
\bbal
||v||_{L^\infty_t(\B^{-1+\frac3p}_{p,1})}&\leq \eta \quad\mbox{for}\quad t\leq \Gamma+\epsilon<T^*,
\end{align*}
which is contradiction with the definition of $\Gamma$.

Thus, we can conclude $\Gamma=T^*$ and
\bbal
||v||_{L^\infty_t(\B^{-1+\frac3p}_{p,1})}&\leq C<\infty \quad\mbox{for all}\quad t\in(0,T^*),
\end{align*}
which implies that $T^*=+\infty$.

\section*{Acknowledgments} J. Li was supported by NSFC (No.11801090).

\end{document}